\documentclass[11pt]{article}
\usepackage{amsfonts,amssymb,amsmath,amsthm}
\usepackage{enumitem}
\usepackage{tikz}

\usepackage{color} 
   \definecolor{cites}{rgb}{0.50 , 0.00 , 0.00}  
   \definecolor{urls} {rgb}{0.00 , 0.00 , 0.50}  
   \definecolor{links}{rgb}{0.00 , 0.00 , 0.50}   
\usepackage[
      colorlinks=true,   
      citecolor=cites,   
      urlcolor=urls,     
      linkcolor=links,   
      pdftitle={Minimal families of limit operators},  
      pdfauthor={Marko Lindner}, 
      pdfpagemode=UseOutlines, 
      pdfstartview=FitH,       
      bookmarksopen=false      
   ]{hyperref}

\newcommand\eps\varepsilon
\newcommand\ph\varphi
\newcommand\spec{{\rm spec}\,}  
\newcommand\speps{{\rm spec}_\eps}
\newcommand\spess{{\rm spec}_{\rm ess}\,}

\newcommand\sppt{{\rm spec}_{\rm pt}^\infty\,}
\newcommand\Lim{{\rm Lim}}

\newcommand\diam{{\rm diam\,}}

\newcommand\clos{{\rm clos}}
\newcommand\ri{{\rm i}}
\newcommand\re{{\rm e}}

\newcommand\im{{\rm im}}
\newcommand\orb{{\rm orb}}
\newcommand\PE{{\rm\Psi E}}
\newcommand\supp{{\rm supp}}
\newcommand\Op{{\rm Op}}

\newcommand\C{{\mathbb C}}
\newcommand\R{{\mathbb R}}

\newcommand\T{{\mathbb T}}
\newcommand\Z{{\mathbb Z}}
\newcommand\N{{\mathbb N}}

\newcommand\M{{\mathcal M}}

\newcommand\cC{{\mathcal C}}
\newcommand\W{{\mathcal W}}
\renewcommand\S{{\Sigma}}

\newtheorem{theorem}{Theorem}[section]
\newtheorem{lemma}[theorem]{Lemma}

\newtheorem{proposition}[theorem]{Proposition}
\newtheorem{definition}[theorem]{Definition}

\newenvironment{example}
 {\par\noindent\refstepcounter{theorem}{\bf Example \thetheorem}\ }
 {\raisebox{1mm}{\framebox{}}\pagebreak[2]}

\newenvironment{Proof}{%
  \begin{proof}[{\bf Proof}]%
}{
   \end{proof}
}

\newtheorem{remarkx}[theorem]{Remark}
\newenvironment{remark}
  {\pushQED{\qed}\remarkx\normalfont}
  {\popQED\endremarkx}

\numberwithin{figure}{section}  

\newcounter{abccounter}
\newenvironment{abc}
  {\begin{list}{\bf \alph{abccounter})}{\usecounter{abccounter}\itemsep-1mm\topsep0mm\parskip1mm}}
  {\end{list}}

\makeatletter
\let\@fnsymbol\@arabic
\makeatother

\begin{document}
\title{\bf Minimal Families of Limit Operators}
\author{
{\sc Marko Lindner}\footnote{Techn. Univ. Hamburg (TUHH), Institut Mathematik, D-21073 Hamburg, Germany, \url{lindner@tuhh.de}}
}

\date{\today}
\maketitle
\begin{quote}
\renewcommand{\baselinestretch}{1.0}
\footnotesize {\sc Abstract.}
We study two abstract scenarios, where an operator family has a certain minimality property.
In both scenarios, it is shown that norm, spectrum and resolvent are the same for all family members. Both abstract settings are illustrated by practically relevant examples, including discrete Schrödinger operators with periodic, quasiperiodic, almost-periodic, Sturmian and pseudo-ergodic potential. The main tool is the method of limit operators, known from studies of Fredholm operators and convergence of projection methods. We close by connecting this tool to the  study of subwords of the operator potential.
\end{quote}

\section{Introduction}
We look at bounded linear operators $A$ on sequence spaces with local action in the sense that the corresponding bi-infinite matrix $(a_{ij})$ is a band matrix. We will often identify the operator and the matrix. We demonstrate everything in the simple situation of $\ell^2(\Z)$ but show in Section~\ref{sec:extend} how to extend the setting to, e.g., vector-valued $\ell^p(\Z^d)$.

A tool that has been established for the study of Fredholm properties, including the essential spectrum, of $A$ is the method of limit operators. Roughly speaking, a (possibly large) family, $\Lim(A)$, of operators, so-called limit operators, is capturing the behaviour of $A$ (or likewise, the asymptotics of the matrix $(a_{ij})$) at infinity. 

In applications ranging from periodic to aperiodic operators \cite{LenzPhD,DamanikLenz2004,Damanik2005,Stollmann}, it turns out that all members of the family $M:=\Lim(A)$ are not only limit operators of $A$ but also of each other, including themselves. 
Moreover, they all share the same spectral quantities like the spectrum itself, the norm, resolvent and pseudospectrum. This observation is not new for the particular classes, although sometimes a bit tedious to derive, but we think it deserves a systematic study in a setting that is based on properties that all classes have in common and that is still strong enough to derive spectral results.
The paper is clearly in the spirit of \cite{Beckus.2017,Beckus.2018} but maybe a bit more concrete, constructive in parts and therefore with more specific results.

By one further step of generalization, we also include so-called pseudo-ergodic operators, first studied by Davies in order to study the spectral theory of random operators while largely eliminating stochastic details.

\section{The tools}
Let us write our band operator $A$ on $\ell^2:=\ell^2(\Z)$ as a finite sum 
\begin{equation} \label{eq:band}
A\ =\ \sum_{j=-w}^w M_{b^{(j)}} S^j
\end{equation}
of products of multiplication operators $(M_{b^{(j)}}x)_n=b^{(j)}_nx_n$ with $b^{(j)}\in\ell^\infty:=\ell^\infty(\Z)$ and powers of the shift operator, $(S x)_n = x_{n - 1}$.
Each coefficient $b^{(j)}$ of $S^j$ represents the $j$-th diagonal of the bi-infinite band matrix behind $A$. By the boundedness of the coefficients, the operator $A$ acts as a bounded linear operator on $\ell^2$; we write $A\in L(\ell^2)$.

Recall that a bounded linear operator $A$ on a Banach space $X$ is a {\sl Fredholm operator}
if its coset, $A+K(X)$, modulo compact operators $K(X)$ is invertible in the so-called Calkin algebra $L(X)/K(X)$. This holds if and only if the nullspace of $A$ has finite dimension and the range of $A$ has finite codimension in $X$.

Since the coset $A+K(X)$ cannot be affected by changing finitely many matrix entries, its study takes place ``at infinity''. This is where limit operators \cite{RaRoSi1998,RaRoSiBook,Li:Book} come in:

\begin{definition}
For a band operator $A$, we look at all its translates $S^{-k}AS^k$ with $k\in\Z$ and speak of a {\sl limit operator}, $A_h$, if, for a particular sequence $h=(h_n)$ in $\Z$ with $|h_n|\to\infty$, the corresponding sequence of translates, $S^{-h_n}AS^{h_n}$, converges strongly to $A_h$.

Here we say that a sequence $A_n$ converges strongly to $A$ if $A_nx\to Ax$ for all $x\in\ell^2$.
\end{definition}

From the matrix perspective, this can be visualized as follows:
$B=(b_{ij})_{i,j\in\Z}$ is the limit operator of $A=(a_{ij})_{i,j\in\Z}$ with respect to the sequence $h=(h_n)$ in $\Z$ if, for all $i,j\in\Z$,
\begin{equation}\label{eq:limop}
a_{i+h_n,j+h_n}\to b_{ij}\quad\textrm{as}\quad n\to\infty.
\end{equation}
We write $A_h$ instead of $B$. Here is the announced connection to Fredholm operators.

\begin{lemma} \label{lem:Fredholm}
For a band operator $A$ on $\ell^2$, it holds that the following are equivalent:\\[-7mm]
\begin{enumerate}[label=(\alph*)]\itemsep-1mm
\item $A$ is a Fredholm operator on $\ell^2$,
\item  all limit operators of $A$ are invertible on $\ell^2$ \cite{RaRoSi1998,LiSei:BigQuest},
\item all limit operators of $A$ are injective on
$\ell^\infty$ \cite{CWLi2008:FC,CWLi2008:Memoir}.
\end{enumerate}
\end{lemma}
Not only because of this result, it is useful to denote the set of all limit operators of $A$ by $\Lim(A)$. Sometimes it is also important to distinguish with respect to the direction: If $h=(h_n)$ with $h_n\to\pm\infty$, we write $A_h\in\Lim_\pm(A)$, respectively.

As a consequence of Lemma \ref{lem:Fredholm},
\begin{equation} \label{eq:specLim}
\spess A\ =\ \bigcup_{A_h\in\Lim(A)} \spec A_h\ =\ \bigcup_{A_h\in\Lim(A)} \sppt A_h.
\end{equation}
One important spectral quantity that we use a lot ``under the hood'' is the so-called {\sl lower norm} of an operator $A$. By this, we mean
\[
\nu(A):=\inf_{\|x\|=1}\|Ax\|.
\]
Note that the name ``lower norm'' is used, as in \cite{RaRoSiBook,Li:Book}, meaning a counterpart of the operator norm. (It is not a norm!) $\nu(A)$ turns out to be a fairly accessible quantity to study $\|A^{-1}\|$. Indeed,
\begin{equation} \label{eq:invnu}
\frac 1{\|A^{-1}\|}\ =\ \min\big(\,\nu(A),\nu(A^*)\,\big).
\end{equation}
We also use that $\nu(A)$ can be conveniently approximated / localized via 
\begin{equation} \label{eq:localnu}
  \nu(A)=\inf_{n\to\infty}\nu_n(A), 
\end{equation}
where $\nu_n(A)=\inf\{\|Ax\|: \|x\|=1, \diam(\supp(x))\le n\}$. For the proof see \cite[Prop. 6]{LiSei:BigQuest} (and \cite[Prop. 3.4]{HagLiSei} for the corresponding result about the norm).

In this way we study resolvent norms, $\|(A-\lambda I)^{-1}\|$, and their superlevel sets, the so-called {\sl $\eps$-pseudospectra} $\speps A=\{\lambda\in\C:\|(A-\lambda I)^{-1}\|>\frac 1\eps\}$, see \cite{TrefEmb}.

\section{Minimal sets of limit operators}
A limit operator of a limit operator of $A$ is a limit operator of $A$, see \cite[Cor. 3.97]{Li:Book}. Precisely,

\begin{lemma} \label{lem:limlim}
For all band operators $A, B$ it holds that
\[
B\in \Lim(A)\quad\implies\quad \Lim(B)\subset \Lim(A).
\]
\end{lemma}
\begin{Proof}
First, $B\in\Lim(A)$ implies $S^{-k}BS^k\in\Lim(A)$ for all $k\in\Z$.
Now let $k=k_1, k_2, ...$, pass to the appropriate limit and note that $\Lim(A)$ is closed in the corresponding sense. For details, see \cite[Cor. 3.97]{Li:Book} or \cite[Lem 3.3]{CWLi2008:FC}.
\end{Proof}
\begin{remark} \label{rem:trans}
{\bf a) }The limit operators of $B$ are even in the same local part, $\Lim_\pm(A)$, of $\Lim(A)$ as $B$ is.
This also generalizes to $\Z^d$, where one has many more directions than left and right.

{\bf b) } Lemma \ref{lem:limlim} shows that the relation
\[
B\prec A\quad :\iff\quad B\in\Lim(A)
\]
is transitive. By \eqref{eq:specLim}, $B\prec A$ implies $\spec B\subset\spess A\subset\spec A$, so that one is also interested in the symmetric (and transitive, but not reflexive) relation
\[
A\sim B\quad :\iff\quad A\prec B\ \ \text{and}\ \ B\prec A,
\]
which yields equality of (essential) spectra and other spectral quantities. We follow this line of thought (on the level of sets $\Lim(A)$ rather than operators $A$) and arrive at operator sets $M$ where any two elements $A,B\in M$ have $A\sim B$ and hence equal spectra.
\end{remark}

\subsection{Two concepts of minimality} \label{sec:31}
Now fix a band operator $A$  and look at the family $\M$ of all sets $M=\Lim(B)$ with $B\in M_0:=\Lim(A)$. By Lemma \ref{lem:limlim}, all these $M$ are subsets of $M_0$ and, starting from $M_0=\Lim(A)$, form a partially ordered set with respect to reversed set inclusion "$\supset$".
By Zorn's lemma (the conditions have been verified in \cite[Prop 3.7]{CWLi2008:FC}), there is a maximal element $M=\Lim(B)\ne\varnothing$ in this poset $(\M,\supset)$ -- let us call it minimal here since the sets get smaller, not bigger. 

As a consequence, for this particular operator set $M$,
\begin{equation} \label{eq:min1}
\forall C\in M:\quad \Lim(C)=M. \tag{M1}
\end{equation}
Our poset construction shows that every set $\Lim(A)$ has a subset $M$ with \eqref{eq:min1}. In some applications, $\Lim(A)$ itself is a set $M$ with \eqref{eq:min1}. This is the case, for example, if just one diagonal of $A$ is varying, and that sequence is\\[-7mm]
\begin{itemize} \itemsep-1mm
\item periodic,
\item almost-periodic, incl.~quasiperiodic, or
\item Sturmian.
\end{itemize}
In other situations, e.g.~a pseudo-ergodic diagonal (we will explain all these classes in Section~\ref{sec:classes} below), we have the following generalization of \eqref{eq:min1}: There is an operator set $M$ such that
\begin{equation} \label{eq:min2}
\forall C\in M:\quad \Lim(C)\supset M. \tag{M2}
\end{equation}
Let us call an operator set $M$ with property \eqref{eq:min2} a {\sl minimal family of limit operators}. Note that \eqref{eq:min1} is a special case of \eqref{eq:min2}.
\medskip

\noindent
\begin{minipage}{0.75\textwidth}
In our derivation of \eqref{eq:min1}, the word ``minimal'' refers to the role of $M$ in $(\M,\supset)$. For completeness, here is the dynamical systems perspective on minimality (also see \cite{Beckus.2017,Beckus.2018}):
For a band operator $B$, call the set of all its translates, that is
\begin{equation} \label{eq:orbB}
\orb(B)\ :=\ \{S^{-n}BS^n:n\in\Z\},
\end{equation}
the {\sl orbit} of $B$. For a set $M$ with \eqref{eq:min1}, resp.~\eqref{eq:min2}, the orbit of every $B\in M$ is dense in $M$, resp.~$M'$ (defined in Proposition \ref{prop:min} b): There are no smaller orbits. 
\end{minipage}
\begin{minipage}{0.25\textwidth}
\begin{tikzpicture}
\draw [gray!30, fill] plot [smooth cycle] coordinates {(0.1,0) (1,0.9) (2,0.8) (2.5,0) (1.1,-0.9) };
\draw[red, thick] (0.4,0) circle (0.3 mm);
\draw (0.3,1) node[below] {$M$};
\draw (3,0.9) node[below] {\eqref{eq:min1}};
\draw [red] plot [smooth cycle] coordinates {
(0.2,0) (1,0.8) (2,0.6) (1.0,-0.8) 
(0.4,0) (1.2,0.7) (2.1,0.5) (1.5,0) (1,-0.3) 
(0.6,-0.1) (1.5,0.6) (2.4,0.3) (1.7,-0.2) (1.2,-0.5)
(0.8,-0.3) (1.6,0.3) (2.3,0) (1.8,-0.5) (1.2,-0.7)};
\draw (-0.5,0);
\end{tikzpicture}
\begin{tikzpicture}
\draw [gray!30, fill] plot [smooth cycle] coordinates {(0.1,0) (1,0.9) (2,0.8) (2.5,0) (1.1,-0.9) };
\draw[red, thick] (0.4,0) circle (0.3 mm);
\draw (0.3,1) node[below] {$M$};
\draw (3,0.9) node[below] {\eqref{eq:min2}};
\draw [red] plot [smooth cycle] coordinates {
(0.2,0) (1,0.8) (2,0.6) (1.0,-0.8) 
(0.4,0) (1.2,0.7) (2.1,0.5) (1.5,0) (1,-0.3) 
(0.6,-0.1) (1.5,0.6) (2.4,0.3) (1.7,-0.2) (1.2,-0.5)
(0.8,-0.3) (1.6,0.3) (2.3,0) (1.8,-0.5) (1.2,-0.7)
(0.6,-0.3) (1.6,0.5) (2.9,-0.1) (2,-1) (1.3,-0.8) 
(1.3,-0.6) (2.2,0) (2.6,-0.2) (1.8,-0.8)};
\draw (-0.5,0);
\end{tikzpicture}

\end{minipage}

Before we prove results about minimal families (and their members) with property \eqref{eq:min1} or \eqref{eq:min2}, let us introduce two more notions that are closely related to \eqref{eq:min1} and \eqref{eq:min2}.

\subsection{Recurrent and self-similar operators} \label{sec:32}
In accordance with \cite{CWLi2008:FC}, and generalizing \cite{Muh1972}, we call an operator $A$ {\sl recurrent} if, for every $B\in\Lim(A)$, one has $\Lim(B)=\Lim(A)$, i.e.~$M:=\Lim(A)$ has property (M1). Moreover, we say that an operator $A$ is {\sl self-similar} if $A\in \Lim(A)$, i.e.~if $A\sim A$.
Some simple interrelations:
\begin{itemize} \itemsep0mm
\item Every limit operator of a recurrent operator clearly is self-similar.

\item Every $B\in M$ with property (M1) is both recurrent and self-similar.

\item Every $B\in M$ with property (M2) is self-similar but might not be recurrent. (e.g.~$B$ pseudo-ergodic, $C\in \Lim(B)$ constant, so that $\Lim(C)\not\ni B$)

\item In particular, self-similar operators need not be recurrent. 

\item But also recurrent operators need not be self-similar.
(e.g.~periodic plus compact)
\end{itemize}

In the simplest case, a minimal set $M$ with \eqref{eq:min1} is a singleton with one translation invariant operator $C$. Here $C$ is called {\sl translation invariant} if $\orb(C)=\{C\}$; from the matrix perspective, this means that all diagonals of $C$ are constant. But $M_0=\Lim(A)$ need not even contain such simple operators and still, a minimal set $M$ exists, by \cite[Prop 3.7]{CWLi2008:FC}.

\begin{lemma}
A minimal set $M$ with \eqref{eq:min1} is a singleton, $M=\{C\}$, if and only if it contains a translation invariant element.
\end{lemma}
\begin{Proof}
If $M$ has just one element, $C$, then $C$ is translation invariant: Otherwise, its translation, $D:=S^{-1}CS^{1}$, were different from $C$ but also in $M=\Lim(B)$. Hence, $M$ had a second element.

If $M$ has more than one element then none is translation invariant: A translation invariant element $C$ would form a ``more minimal'' subset $\tilde M:=\{C\}=\orb(C)=\Lim(C)\subsetneq M$.
\end{Proof}
The next two lemmas follow immediately from Lemma \ref{lem:Fredholm}.
\begin{lemma} \label{lem:1-4}
For a  set $M$ with property \eqref{eq:min1}, the following conditions are equivalent:
\begin{minipage}{6mm}
~
\end{minipage}
\begin{minipage}{0.7\textwidth}
\begin{enumerate}[label=$(\roman*)$\ ] \itemsep-1mm
\item one $C\in M$ is a Fredholm operator on $\ell^2$,
\item one $C\in M$ is invertible on $\ell^2$,
\item all $C\in M$ are invertible on $\ell^2$,
\item all $C\in M$ are injective on $\ell^\infty$.
\end{enumerate}
\end{minipage}
\end{lemma}
%
$(i)\Leftrightarrow(ii)$ holds on the level of individual self-similar operators, by Lemma \ref{lem:Fredholm} and \eqref{eq:specLim}:
\begin{lemma}
If $A$ is self-similar, e.g.~if $A\in M$ with \eqref{eq:min2}, then
\begin{abc}
\item $A$ is invertible if and only if $A$ is Fredholm,
\item $\spec A=\spess A$.
\end{abc}
\end{lemma}
Note that, similar to $(iii)\Leftrightarrow(iv)$ in Lemma \ref{lem:1-4}, equivalence between invertibility on $\ell^2$ and injectivity on $\ell^\infty$ can be shown -- even for individual operators and under much weaker conditions. However, in that weaker setting, also injectivity of $A^*$ is required (which is of course redundant if $A$ is self-adjoint):

\begin{lemma} \label{lem:inj}
{\bf a) } For a band operator $A$ with closed range (e.g.~a Fredholm operator), one has:
\[
\text{$A$ is invertible on $\ell^2 \qquad\iff\qquad A$ and $A^*$ are injective on $\ell^\infty$.}
\]
{\bf b) } Similarly, the following equivalence holds for the half-line compression $A_+$ of $A$:
\[
\text{$A_+$ is invertible on $\ell^2(\Z_+) \qquad\iff\qquad A_+$ and $A_+^*$ are injective on $\ell^\infty(\Z_+)$.}
\]
Here $A_+$ denotes the half-line compression of $A$, that is $PAP$, as an operator $\ell^2(\Z_+)\to\ell^2(\Z_+)$, where $P:\ell^2\to\ell^2$ is the operator of multiplication by the characteristic function of $\Z_+$.
\end{lemma}
\begin{Proof}
{\bf a) } $\boxed\Rightarrow$ If $A$ is invertible on $\ell^2$, then so is its adjoint, $A^*$. As band operators, $A$ and $A^*$ act as bounded operators on every $\ell^p$ with $p\in[1,\infty]$, where, by \cite{RaRoSiBook,Li:Habil}, their invertibility and spectrum are independent of $p\in[1,\infty]$. So both are invertible and hence injective on $\ell^\infty$.

$\boxed\Leftarrow$ If $A$ and $A^*$ are injective on $\ell^\infty$, then they are also injective on $\ell^2\subset\ell^\infty$. Moreover, the range of $A$ in $\ell^2$ is closed, by assumption. 

{\bf b) } 
The proof of b) works exactly like a). The range of $A_+$ in $\ell^2(\Z_+)$ is closed, too, since $\im(A_+)=\im(PAP)=\im(P)\cap\im(AP)$ with $\im(AP)=\im(PA)+\im(\chi_{\{-w,...,-1\}}\cdot)$, where $w$ is the band-width of $A$ and $\im(PA)=\im(P)\cap\im(A)$.
\end{Proof}

Note that statements a) and b) of Lemma \ref{lem:inj} also hold with injectivity on any $\ell^p$ with $p\ge 2$. But of course, technically, $p=\infty$ is the simplest case to check. ($A_+x=0$ yields a componentwise recurrence for a vector $x$ in the kernel; now test for boundedness).

\begin{proposition} \label{prop:min}
Let an operator set $M$ be subject to \eqref{eq:min2}.
Then
\begin{abc}
\item all elements of $M$ are limit operators of each other -- including themselves,
\item all sets $\Lim(C)$ in \eqref{eq:min2} are the same; denote that set by $M'$ ($=M$ if even \eqref{eq:min1} applies),
\item all elements of $M$ have the same norm, lower norm, spectrum (that is entirely essential spectrum), resolvent norms and pseudospectra,
\item for all $C\in M$, one has $\Lim_-(C)=\Lim(C)=\Lim_+(C)$.
\item for all $C\in M$, one has $\spess C_-=\spess C_+=\spess C=\spec C$, where $C_\pm$ denotes the compression of $C$ to the positive resp.~negative half-line, $\ell^2(\Z_\pm)$.
\end{abc}
\end{proposition}
\begin{Proof}\noindent
\begin{abc}
\item If $A, B\in M$ then, by \eqref{eq:min2}, $A\in M\subset\Lim(B)$ and $B\in M\subset\Lim(A)$.
\item If $A, B\in M$ then, by a) and Lemma \ref{lem:limlim}, $\Lim(A)\subset\Lim(B)$ and $\Lim(B)\subset\Lim(A)$.

\item $B\in\Lim(A)$ implies $\|B\|\le\|A\|$ by \cite[Prop 1.a]{RaRoSi1998}, $\nu(B)\ge\nu(A)$, by \cite[Prop 3.9]{HagLiSei} and $\spec B\subset\spess A\subset \spec A$ by \eqref{eq:specLim}. By a), all $A,B\in M$ are limit operators of each other.

Moreover,  if $A,B\in M$, so that $B=A_h$ for some $h$, by a), then $B-\lambda I=A_h-\lambda I=(A-\lambda I)_h$ for all $\lambda\in\C$, so that $A-\lambda I$ and $B-\lambda I$ are limit operators of each other. Hence, also $\nu(A-\lambda I)$ and $\nu(B-\lambda I)$ coincide -- as well as their level sets.

\item Let $C\in M$ and $B\in \Lim_+(C)$. Then, by Remark \ref{rem:trans} a), and recalling $M'$ from b), 
\[
M'=\Lim(B)\subset \Lim_+(C)\subset \Lim(C)=M',
\]
so that both ``$\subset$'' are ``$=$'', whence $\Lim_+(C)=\Lim(C)$. It is the same with $\Lim_-(C)$.

\item One defines limit operators of one-sided infinite matrices $C_\pm=(C_{ij})_{i,j\in\Z_\pm}$ (in the direction where the matrix is infinite) by the same definition and gets $\Lim(C_\pm)=\Lim_\pm(C)$. By d), both sets are equal to $\Lim(C)$, and by \eqref{eq:specLim} and its one-sided analogue, it follows $\spess C_\pm=\spess C$. Equality with $\spec C$ is already in c).
\qedhere 
\end{abc}
\end{Proof}

Let $S_M\subset \mathbb C$ denote the (essential) spectrum of one (and all) operator(s) $C\in M$.
The sets $\spec C_+$ and $S_M$ differ by eigenvalues of $C_+$ or its adjoint -- so-called {\sl Dirichlet eigenvalues}. Note that this difference set, and hence $\spec C_+$ (or $\spec C_-$), can be different or the same for all $C\in M$ -- depending on the family $M$, e.g.
\begin{itemize}
\item For a family $M$ of periodic discrete Schrödinger operators (i.e.~the orbit of one $q$-periodic discrete Schrödinger operator containing $q$ elements), there is at most one Dirichlet eigenvalue per gap between the (at most) $q$ intervals ({\sl spectral bands}) of $S_M\subset\R$. But these eigenvalues typically differ for every $C\in M$.

\item For tridiagonal pseudoergodic operators, the difference between $\spec C_+$ and $S_M$ is the same for all $C\in M$, see \cite{CWLi:Coburn} for a proof and the explicit computation of that difference.
\end{itemize}
In particular, it is in general impossible to make uniform statements for the whole family $M$ about the invertibility of $C_\pm$ or about more involved constructions like the finite section method.

\subsection{Concrete operator classes and examples} \label{sec:classes}
Let us look at classes of recurrent operators $A$, so that already the set $M:=\Lim(A)$ is minimal in the sense  
\eqref{eq:min1}. 
\medskip

{\bf Periodic and almost-periodic operators.\ }
A band operator $A$ is called a {\sl periodic}, resp.~{\sl almost-periodic}, {\sl operator} if $\orb(A)$ is finite, resp.~relatively compact in the uniform (meaning operator norm) topology.

One can show, see \cite{Cordu} or \cite[Lem. 5.43]{Li:Habil}, that $A$ is a periodic, resp.~almost-periodic, operator if and only if every coefficient $b^{(j)}$ in \eqref{eq:band} is periodic, resp. almost-periodic. Here, a sequence $b$ is called {\sl periodic}, resp.~{\sl almost-periodic}, if its {\sl orbit}
\begin{equation} \label{eq:orbb}
\orb(b) \ :=\ \{S^n b\ :\ n\in\Z\}
\end{equation}
is finite, resp.~relatively compact in $\ell^\infty(\Z)$.

\begin{itemize}
\item
For a $q$-periodic operator $A$ with $q\in\N$, we have
\[
\Lim(A)\ =\ \{S^{-n}AS^n:n=0,\dots,q-1\}\ =\ \orb(A).
\]
So the set $M:=\Lim(A)$ is minimal in the sense \eqref{eq:min1}, hence \eqref{eq:min2}, so that $M$ is subject to Proposition \ref{prop:min}. Of course, none of the statements of Proposition \ref{prop:min} is a surprise in this setting, as every $C\in M$ is related to $A$ by a unitary similarity transform, $C=S^{-k}AS^k$.

\item
The situation is slightly different for almost-periodic operators $A$:
If $B\in\Lim(A)$ then the convergence $S^{-h_n}AS^{h_n}\to B$ is, unlike in the periodic case, not by being eventually constant but, still remarkably, in the operator norm topology!
Indeed, the sequence $(S^{-h_n}AS^{h_n})$ has, by the relative compactness of $\orb(A)$, a uniformly convergent subsequence, clearly with limit $B$.

Put $M:=\Lim(A)$ and let us check \eqref{eq:min1}. If $C\in M$ then $C=A_h$ with a sequence $h=(h_n)$ in $\Z$, where $|h_n|\to\infty$.
As the convergence $S^{-h_n}AS^{h_n}\to C$ is even uniform, we get, 
\[
\|A-S^{h_n}CS^{-h_n}\|\ =\ \|S^{-h_n}AS^{h_n}-C\|\ \to\ 0,
\]
so that $A=C_{-h}\in\Lim(C)$, whence $M=\Lim(A)\subset\Lim(C)$, by Lemma \ref{lem:limlim}.
The opposite inclusion, $\Lim(C)\subset M$ holds by $C\in M$ and Lemma \ref{lem:limlim}.
So indeed, \eqref{eq:min1} holds.

Also in the almost-periodic setting, one can conclude
$\|A\|=\|S^{-h_n}AS^{h_n}\|\to\|B\|$ and $\nu(A)=\nu(S^{-h_n}AS^{h_n})\to \nu(B)$, hence equality, directly from the uniform convergence $S^{-h_n}AS^{h_n}\to B$. The statements for spectra, resolvents, etc.~follow easily also from here.
\end{itemize}

\begin{example}
Prominent examples are so-called {\sl discrete Schrödinger operators}
\begin{equation} \label{eq:Schroed}
A\ =\ S^{-1}+M_{b^{(0)}}+S,
\end{equation}
where $b:=b^{(0)}$ is called the {\sl potential} of $A$. Such operators appear, for example, in condensed matter physics, describing the electrical conductivity of materials with regular (periodic, i.e.~crystal) or less regular structure.

A typical construction of potentials $b$ is as follows. For $n\in\Z$, let
\begin{equation} \label{eq:AM}
b(n)=g(f(n))
\qquad\text{with}\qquad
f(x)=\re^{2\pi\ri\,(\theta + \alpha x)}
\qquad\text{and}\qquad
g(t)=\text{Re}(t).
\end{equation}
If $\alpha$ is rational, $\frac pq$, then $f$, and hence $b$, is $q$-periodic (assuming the fraction $\frac pq$ cannot be reduced).
For irrational $\alpha$, however, $f(n)$ never repeats a value, and $b=g\circ f$ is not periodic.
But, due to the compactness of the unit circle $\T:=f(\R)$ and the continuity of $g$, it is easy to see that every sequence of elements of $\orb(b)$ has a uniformly convergent subsequence.

The discrete Schrödinger operator with potential $b$ from \eqref{eq:AM} (or a multiple of $b$) is the 
well-studied {\sl Almost-Mathieu operator}. The particular construction  $b=g\circ f$ in \eqref{eq:AM} with a continuous $g$ is also referred to as a {\sl quasi-periodic} potential.

Let $A$ be the Almost-Mathieu operator. We know from above that $M:=\Lim(A)$ has property \eqref{eq:min1}. For this particular construction, it is known that $M$ is the set of all Schrödinger operators with potential $b$ of the form \eqref{eq:AM} and $\theta\in [0,1)$ -- regardless of the $\theta$ used in $A$.
\end{example}
\medskip

{\bf Sturmian models. \ } 
This is a class, where the limit operators and the fact \eqref{eq:min1} are well-known \cite{LenzPhD,DamanikLenz2004,Damanik2005,Stollmann} and where it is much more tedious than in the periodic or almost-periodic case to conclude the claims of Proposition \ref{prop:min} by bare hands.

For simplicity, again look at a band operator with just one varying diagonal, e.g.~a discrete Schrödinger operator, and copy the construction from \eqref{eq:AM} but, in contrast, replace $g$ by a discontinuous function, e.g.~the characteristic function $g=\chi_I$ of the circular interval $I=f([1-\beta,1))$. 

The model then has three parameters: the step size\footnote{Again, a rational $\alpha$ leads to periodicity, whence $\alpha$ is typically chosen irrational.} $\alpha$ in $f$, the interval length $\beta$ in $g$ and the initial position $\theta$ in $f$. Another degree of freedom is to replace the interval $I=f([1-\beta,1))$ by $J=f((1-\beta,1])$. Let us denote the corresponding operator $A$ by $A_{\alpha,\beta}^{\theta,I}$ resp.~$A_{\alpha,\beta}^{\theta,J}$.

All these configurations produce so-called {\sl Sturmian} sequences, and in all cases,
\[
\Lim(A_{\alpha,\beta}^{\theta,I})\ 
=\ \Lim(A_{\alpha,\beta}^{\theta,J})\ 
=\ \left\{\ A_{\alpha,\beta}^{\phi,I},\ A_{\alpha,\beta}^{\phi,J}\ \ :\ \ \phi\in [0,1)\ \right\}\ =:\ M,
\]
independently of $\theta$. So obviously, \eqref{eq:min1} holds, whence Proposition \ref{prop:min} is in force.

The standard model, the so-called {\sl Fibonacci Hamiltonian}, has $\theta=0$, $g=\chi_I$, and
\begin{equation} \label{eq:FiboHam}
\alpha\ =\ \beta\ =\ \frac{\sqrt 5-1}2\ =\ \frac 1{1+\frac 1{1+\frac 1{1+\cdots}}}.
\end{equation}

{\bf Pseudo-ergodic models.\ }
Now we turn to a case where \eqref{eq:min1} does not apply -- but \eqref{eq:min2}.

For simplicity, again look at a band operator $A$ with just one varying diagonal $b$, e.g.~a discrete Schrödinger operator. Denote this $A$ by $A_b$. Also note our construction in Section \ref{sec:extend} below how to extend this to finitely many varying diagonals.

For the entries of $b$, fix a compact set $\S\subset\C$ -- the so-called {\sl alphabet}.
We say that both $b\in\S^\Z$ and $A_b$ are {\sl pseudo-ergodic}, and write $A_b\in\PE$, if every $A_c$ with $c\in\S^\Z$ is a limit operator of $A_b$. Abbreviating the set $\{A_c : c\in \S^\Z\}$ by $A_{\S^\Z}$, we get that
\begin{equation} \label{eq:PE}
\forall A\in\PE:\quad \Lim(A)\ =\ A_{\S^\Z}\ \supset\ \PE,
\end{equation}
so that $M:=\PE$ is subject to \eqref{eq:min2} and Proposition \ref{prop:min} is in force.
In particular, $M'$ is $A_{\S^\Z}$, the set of all operators $A_c$ that are possible with a sequence $c$ over the alphabet $\S$.

It is also possible to construct minimal families $M$ with \eqref{eq:min2}, where $M\subsetneq M'\subsetneq A_{\S^\Z}$. To this end, assume $|\S|\ge 2$, say $0,1\in\S$, and limit the basic building blocks of the bi-infinite word $b$ to a set $T$ of finite words, 
e.g.~$T=\{00,1\}$, meaning that we only consider $b\in\S^\Z$ that arise from bi-infinite concatenations of ``allowed'' words $t\in T$. 
\pagebreak

Denote the corresponding set of operators $A_b$ by $\Op(T)$. In particular, $\Op(\S)=A_{\S^\Z}$. If $\Op(T)\subsetneq\Op(\S)$ then, by limiting ourselves to pseudo-ergodic operators in $\Op(T)$, \eqref{eq:PE} generalizes as follows:
\begin{equation} \label{eq:PET}
\forall A\in\underbrace{\PE\cap \Op(T)}_M:\quad A_{\S^\Z}\ \supsetneq\ \Lim(A)\ =\ \underbrace{A_{\S^\Z}\cap \Op(T)}_{M'}\ \supsetneq\ \underbrace{\PE\cap \Op(T)}_M\,.
\end{equation}
Pseudo-ergodicity was introduced by Davies \cite{Davies2001:PseudoErg} to study spectral properties of random operators while eliminating probabilistic arguments. Indeed, if the $b(n)$ are random variables then, for a large class of probability distributions, in particular if each $b(n)$ is iid with range $\S$, then $A_b\in\PE$ almost surely.

\section{Subwords a.k.a.~factors}
The concept of pseudo-ergodicity (and some others) is much more user-friendly when defined in terms of subwords rather than limit operators: \ $b\in\S^\Z$ is pseudo-ergodic if and only if
\begin{center}
{\em
every finite word over $\S$ can be found, up to arbitrary precision, as a subword of $b$.
}
\end{center}

If $\S$ is a {\sl discrete alphabet}, meaning that
\[
\inf_{s,t\in\S,\ s\ne t}|s-t|\ >\ 0,
\]
then the ``up to arbitrary precision'' bit can even be dropped in the previous sentence. Of course:
\begin{lemma} \label{lem:discrete}
{\bf a) } If $\S$ is discrete then a sequence in $\S$ is convergent iff it is eventually constant.\\
\indent
{\bf b) } A bounded alphabet $\S\subset\C$ is discrete if and only if it is finite. (Bolzano-Weierstrass)
\end{lemma}
We should now define the notion of a {\sl subword} and some related notations properly:\\[-7mm]
\begin{itemize} \itemsep-1mm
\item Let $\S^* := \cup_{n=0}^\infty\ \S^n$ denote the set of all finite words over $\S$.
\item Let $[a_1\dots a_n]$ denote the word $w\in\S^n$ with $w(k)=a_k\in\S$ for $k=1,...,n\in\N$.
\item Let $|w|$ be the length of the word $w\in\S^*$; so $|[a_1\dots a_n]|=n$.
\item For $w\in\S^*$, $b\in\S^\Z$ and $\eps>0$, put
\\[-7mm]
\begin{itemize} \itemsep-1mm
\item[$\circ\ $] $pos(w,b):=\{k\in \Z: [b(k)\dots b(k+|w|-1)]=w\}$,
\item[$\circ\ $] $pos_\eps(w,b):=\{k\in \Z: \|\,[b(k)\dots b(k+|w|-1)]-w\,\|_\infty<\eps\}$,
\item[$\circ\ $] $\#(w,b):=|pos(w,b)|$,\ \ meaning the number of elements,
\item[$\circ\ $] $\#_\eps(w,b):=|pos_\eps(w,b)|$,
\item[$\circ\ $] $\W(b) \ :=\ \{w\in\S^*:\#(w,b)\ge 1\}$, the set of all finite subwords of $b$, and
\item[$\circ\ $]  $\W_n(b):=\{w\in\W(b):|w|=n\}$.
\end{itemize}
\end{itemize}

The two concepts of limit operators and subwords (a.k.a.~factors) go hand in hand very nicely in the case of pseudo-ergodicity. Also: both concepts yield neat and clean sufficient conditions for inequalities and equalities between norms, lower norms, spectra, pseudospectra, etc.~-- see Proposition \ref{prop:min} above and Proposition \ref{lem:subw-nu} below. So let us also have a look at how they are related to each other in general.

\subsection{Subwords and spectra}
First, here is how subwords indicate spectral inclusions: Let $A_b$, again, be a band operator with one distinguished diagonal carrying the sequence $b\in\S^\Z$ and all other diagonals constant. (And~note our construction in Section \ref{sec:extend} on how to extend this to finitely many varying diagonals.) 
\begin{proposition} \label{lem:subw-nu}
If $b,c\in\S^\Z$ with $\W(b)\subset\W(c)$ then 
\begin{abc}
\item $\nu(A_b)\ge\nu(A_c)$,
\item $\spec A_b\subset\spec A_c$ and
\item $\speps A_b\subset\speps A_c$ for all $\eps>0$.
\end{abc}
\end{proposition}

\begin{Proof}
Our attention was drawn to observation a) by \cite{Fab}. It will also be part of \cite{aperSchr}. Let's give the simple proof here, for the reader's convenience:
\begin{abc}
\item Take any $x\in\ell^2$ with $\|x\|=1$ and finite support, say of size $n\in\N$. 
As $\W(b) \subset \W(c)$, there exists $k \in \Z$ such that $\|A_b x\| = \|A_c S^k x \|$,
which implies $\nu_n(A_b) \geq \nu_n(A_c)$ for all $n$. Now use \eqref{eq:localnu} to conclude $\nu(A_b)\ge\nu(A_c)$.

\item $\W(b)\subset\W(c)$ and $\lambda\in\C$ imply $\W(b-\lambda)\subset\W(c-\lambda)$. Whether it is the main diagonal that is varying or not, $\nu(A_b-\lambda I)\ge\nu(A_c-\lambda I)$, by part a). By the same arguments, $\nu((A_b-\lambda I)^*)\ge\nu((A_c-\lambda I)^*)$ holds, so that $\|(A_b-\lambda I)^{-1}\|\le\|(A_c-\lambda I)^{-1}\|$, by \eqref{eq:invnu}.

\item This follows immediately from $\|(A_b-\lambda I)^{-1}\|\le\|(A_c-\lambda I)^{-1}\|$.
\qedhere
\end{abc}
\end{Proof}

\subsection{Subwords and limit operators}
Now we connect the concepts of subwords and limit operators, as good as we can. A perfect match is not to be expected, though, since the limit operator relation, $A_b\in\Lim(A_c)$, is connected with infinite repetition (or approximation) of all finite patterns of $b$ in $c$. We begin by making this observation explicit.
\begin{lemma} \label{lem:limopW}
Let $b,c\in\S^\Z$. Then the following are equivalent\\[-7mm]
\begin{enumerate}[label=(\roman*)\ ] \itemsep-1mm
\item $A_b\in \Lim(A_c)$,
\item for every $w\in\W(b)$, it holds that $\#_\eps(w,c)=\infty$ for all $\eps>0$; that is, either\\[-7mm]
\begin{enumerate}[label=(\alph*)\ ]
\item $\#(w,c)=\infty$\ \ or
\item $w\in\clos\big(\W_n(c)\setminus\{w\}\big)$ with $n=|w|$.
\end{enumerate}
\end{enumerate}
If $\S$ is discrete then (b) is impossible, so that (ii) is always (a) then.
\end{lemma}
\begin{Proof}
$\boxed\Downarrow$ Let $A_b\in \Lim(A_c)$, $w\in\W(b)$ and $k\in pos(w,b)$, so that $w=[b(k)\dots b(k+|w|-1)]$.
By $A_b\in\Lim(A_c)$, $A_b=(A_c)_h$ for a sequence $h=(h_n)$ in $\Z$ with $|h_n|\to\infty$. In particular,
\[
w=\lim_{n\to\infty} [c(k+h_n)\dots c(k+|w|-1+h_n)].
\]
Now we have two cases: The sequence $(w_n)$ of words $w_n=[c(k+h_n)\dots c(k+|w|-1+h_n)]\in\S^n$ \\[-7mm]
\begin{enumerate}[label=$(\alph*)$\ ] \itemsep-1mm
\item ... contains $w$ infinitely often (i.e.~it has a constant subsequence). Then $\#(w,c)=\infty$.
\item  ... contains $w$ only finitely often. Then we can remove those finitely many $w$ from $(w_n)$, still have $w_n\to w$ and conclude that $w\in\clos(\W_{|w|}(c)\setminus\{w\})$.
\end{enumerate}

$\boxed\Uparrow$  For $n\in\N$, choose $h_n\in pos_{1/n}([b(-n)\dots b(n)],c)$. By $(ii)$, the choice is infinite; so let $|h_n|>n$. 
Then $|h_n|\to\infty$ and $S^{-h_n}A_cS^{h_n}\to A_b$ strongly, 
whence $A_b=(A_c)_h\in \Lim(A_c)$.
\end{Proof}
\pagebreak

As a consequence, we immediately get that for two members $A_b$ and $A_c$ of a minimal family $M$ with \eqref{eq:min2}, in the case of a discrete alphabet $\S$ for $b$ and $c$, we have $\W(b)=\W(c)$ and that every $w\in\W(b)$ appears infinitely often in $b$ -- in the left and in the right half of $b$.

\begin{remark}
With Lemma \ref{lem:limopW}, we can go and compare the two conditions\\[-7mm]
\begin{enumerate}[label=$(\roman*)$\ ] \itemsep-1mm
\item $A_b\in \Lim(A_c)$,
\item $\W(b)\subset\W(c)$.
\end{enumerate}
Both are sufficient for $\spec A_b\subset\spec A_c$, etc., by Propositions \ref{prop:min} and \ref{lem:subw-nu}. The relation to each other is more shaky. In fact, we cannot expect an immediate connection since $\Lim(A_b)$ determines the essential spectrum of $A_b$, by \eqref{eq:specLim}, while $w\in\W(b)$ with $\#_\eps(w,b)<\infty$ does not.

$(i)\Rightarrow (ii):$ This holds if $\S$ is discrete (even with $\#(w,c)=\infty$). Otherwise, $(i)$ implies at least $\W(b)\subset \clos\,\W(c):=\cup_{n=0}^\infty \clos\, \W_n(c)$.

$(ii)\Rightarrow (i):$ Let $w\in\W(b)\subset\W(c)$.
For $A_b\in\Lim(A_c)$ we even need $\#_\eps(w,c)=\infty$ for all $\eps>0$, 
by Lemma \ref{lem:limopW}.
So this direction needs further strong conditions. Let us assume $\S$ is discrete.
For the conclusion  $\#(w,c)=\infty$, it helps, for example, when $A_b$ or $A_c$ is self-similar:
\[
\begin{array}{ccccc}
\#(w,b)\ge 1 &
\stackrel{A_b\in\Lim(A_b)}\implies &
\#(w,b)=\infty &
\stackrel{\W(b)\subset\W(c)}\implies &
\#(w,c)=\infty,
\\
\#(w,b)\ge 1 &
\stackrel{\W(b)\subset\W(c)}\implies &
\#(w,c)\ge 1 &
\stackrel{A_c\in\Lim(A_c)}\implies &
\#(w,c)=\infty. \qedhere
\end{array}
\]
\end{remark}

\begin{remark} \label{rem:LinRep}
{\bf (Bounded gaps and linear repetitivity)\ }
Damanik and Lenz \cite{DamanikLenz2006} study minimal families $M=\Lim(A_b)$ over a finite (hence discrete) alphabet. We know that then every $A_c\in M$ is self-similar, i.e. every subword $w\in\W(c)$ occurs $\#(w,c)=\infty$ often.
It is shown that \\[-7mm]
\begin{enumerate}[label=$(\roman*)$\ ] \itemsep-1mm
\item $M$ is minimal in the sense that the orbits of all $A_c\in M$ are dense in $M$, if and only if

\item the gap length between any two occurrences of $w$ in $c$ is bounded (the bound dependent on $w$ but not on $A_c\in M$), which obviously follows from

\item the gap between any two occurrences of $w$ in $c$ is bounded by $L|w|$ for all $A_c\in M$ with a constant $L$ independent of $w$ and $c$.
\end{enumerate}
Property $(iii)$ is called {\sl linear repetitivity} of $M$ and is shown to be also necessary for $(i)$ and $(ii)$ if the infinite word $b$ is generated by substitution rules, like the Fibonacci potential \eqref{eq:FiboHam} is recursively generated by application of the two rules $0\mapsto 1,\ 1\mapsto 10$. 

Clearly, minimality $(i)$ here means that $\clos(\orb(A_c))$ equals $M$ for all $A_c\in M$ but cannot be larger than $M$; in our notations, that is (M1) but not its generalization (M2). Indeed:\\[-7mm]
\begin{itemize} \itemsep-1mm
\item boundedness of the gap between two occurrences of $w$ in $c$ implies that $w$ also occurs in every $d$ with $A_d\in\Lim(A_c)$ -- infinitely often and with the same bound on the gaps; limit operators $A_d$ without or with just finitely many $w$ in $d$ are impossible,

\item in particular, pseudo-ergodic cases are not covered: every $w\in\S^*$ occurs infinitely often in every pseudo-ergodic $c\in\S^\Z$ but we have no control over the locations and hence no bound on the gap size. Unbounded gaps lead to limit operators $A_d$ of $A_c$ with $\#(w,d)=0$.
\qedhere
\end{itemize}
\end{remark}
 
\subsection{Factor complexity and operator classes} 
Finally, note that some of our classes of operators $A_b$ and sequences $b\in\S^\Z$ can be characterized purely by the size of the set $\W_n(b)$ -- the so-called {\sl factor complexity} -- that is the number of different subwords of $b$ of length $n$. Here are the basics (see e.g.~\cite{CovenHedlund}):
\begin{itemize} \itemsep0mm
\item The function $\cC_b(n):=|\W_n(b)|$ is only non-trivial for $|\S|>1$.
\item The function $\cC_b(n)$ is monotonic non-decreasing. 
\item Once the function stops growing, $\cC_b(n)=\cC_b(n+1)$, it is constant from there.
\item So if $\cC_b$ is unbounded then $\cC_b(n)\ge n+1$.
\item If $b$ is periodic, say with period $q$, then $\cC_b(n)\le q$, with ``$=$'' if $n\ge q$.
\item In fact, $\cC_b$ is bounded if and only if $b$ is periodic.
(bi-infinite Moore-Hedlund theorem, \cite{CovenHedlund})
\item A one-sided infinite word $a$, e.g.~the left half, $b_-$, or the right half, $b_+$, of a bi-infinite word~$b$,  is Sturmian if and only if $\cC_{a}(n)=n+1$ (minimal unbounded growth).
\item For bi-infinite words, the situation is more complicated, see \cite{CovenHedlund}, and note that $\cC_b(n)=n+1$ still allows bi-infinite words like $b=\dots000111\dots$ or $\dots0001000\dots$
\item For Sturmian words $b$, one has $|\S|=\cC_{b_\pm}(1)=2$, so that the choice $\Sigma=\{0,1\}$ is not as arbitrary as it first seems.
\item $b$ is pseudo-ergodic if and only if $\cC_b(n)=|\S|^n$.
\end{itemize}

\section{Possible directions of extension} \label{sec:extend}
\begin{itemize}
\item We can clearly pass from $\ell^2$ to $\ell^p$, even with $p\in[1,\infty]$.
Instead of the strong convergence of $S^{-h_n}AS^{h_n}$, one then looks at the so-called $\mathcal P$-convergence \cite{Sei:Survey}. In particular, none of our results has anything to do with self-adjointness.

\item Our sequence spaces could also have values in a Banach space $X$. The matrix entries of $A$ (and hence the elements of $\S$) will then have to be operators $X\to X$. Also this setting is covered by the $\mathcal P$-convergence and corresponding results. One should however be aware that then, unlike for $\dim X<\infty$, not for every bounded band operator $A$ and every prescribed sequence $h=(h_n)$ in $\Z$ with $|h_n|\to\infty$, there is a subsequence $g$ of $h$ for which the limit operator $A_g$ exists.

\item For simplicity, we restricted consideration to band operators with just one varying diagonal.
One can also study finitely many, say $m$, varying diagonals:
Therefore write $A=A_b$ with a sequence $b\in T^\Z$, where $T=\Sigma^m$ and $b(k)$ is the vector containing the $m$ interesting entries in the $k$-th column of $A$.

\item In Sections \ref{sec:31} and \ref{sec:32}, we can pass from $\ell^p(\Z)$ to $\ell^p(\Z^d)$. Of course, even 2D quasicrystals, the generalization of Sturmian sequences, the most famous one being the Penrose tiling, are and remain a discipline of its own. Also for simpler classes, one has to adapt the notion of a 2D or 3D subword accordingly, of course.

\item Speaking about Sturmian sequences, there exist generalizations to alphabets $\S$ with $k\ge 2$ elements. The condition is then a factor complexity of $n+k-1$.

\item One can also pass from band operators to uniform limits of band operators -- so-called {\sl band-dominated operators}, \cite{Li:Habil}.
\end{itemize}

\medskip

{\bf Acknowledgements.}
The author thanks Fabian Gabel, Dennis Gallaun, Julian Großmann, Christian Seifert and Riko Ukena for helpful comments and discussions.

\vfill
\noindent {\bf Author's address:}\\
\\
Marko Lindner\hfill \href{mailto:lindner@tuhh.de}{\tt lindner@tuhh.de}\\
Institut Mathematik\\
TU Hamburg (TUHH)\\
D--21073 Hamburg\\
GERMANY
\end{document}